\documentclass{elsart3-1}

 \usepackage{graphicx}
 \usepackage{epsfig}

\usepackage{amssymb,MnSymbol,cite}
\usepackage{mathrsfs}
\usepackage[english,francais]{babel}
\usepackage{caption}
\usepackage{subcaption}

\newtheorem{theorem}{Theorem}[section]
\newtheorem{lemma}[theorem]{Lemma}
\newtheorem{e-proposition}[theorem]{Proposition}

\newtheorem{e-definition}[theorem]{Definition\rm}

\renewcommand{\theequation}{\arabic{equation}}
\def\og{\leavevmode\raise.3ex\hbox{$\scriptscriptstyle\langle\!\langle$~}}
\def\fg{\leavevmode\raise.3ex\hbox{~$\!\scriptscriptstyle\,\rangle\!\rangle$}}

\def\Hom{{\rm Hom}}

\def\Spec{{\rm Spec\,}}

\def\A{{\mathbb A}}
\def\B{{\mathbb B}}
\def\C{{\mathbb C}}

\def\N{{\mathbb N}}

\def\Q{{\mathbb Q}}
\def\R{{\mathbb R}}
\def\Z{{\mathbb Z}}

\def\aarith{{\mathscr A}}
\def\scal{{(\rnt,\cO)}}
\def\scal1{{\hat \aarith}}

\def\cdim{{{\mbox{Dim}_\R}}}
\def\tdim{{{\mbox{dim}_{\rm top}}}}

\def\cE{{\mathcal E}}

\def\cK{{\mathcal K}}

\def\cO{{\mathcal O}}
\def\cP{{\mathcal P}}

\def\cR{{\mathcal R}}

\def\cU{{\mathcal U}}

\def\cZ{{\mathcal Z}}

\def\qqq{\,,\,~\forall}

\newcommand{\ie}{{\it i.e.\/}\ }

\newcommand{\cf}{{\it cf.}}

\def\Hom {{\mbox{Hom}}}

\def\ffp{\mathfrak{p}}

\def\mc{multiplicatively cancellative }
\def\zmax{{\Z_{\rm max}}}

\def\rmax{\R_+^{\rm max}}
\def\rma{\R_{\rm max}}

\def\Se{\frak{ Sets}}

\def\tt{\frak{t}}

\def\sh{\mathfrak{Sh}}

\newcommand{\nil}[1]{}

\def\nt{\N^{\times}}

\def\wnt{{\widehat{\N^{\times}}}}
\def\rnt{{[0,\infty)\rtimes{\N^{\times}}}}

\def\div{{\rm Div}}

\def\vsp{\vspace{.05in}}

\newcommand\blfootnote[1]{%
  \begingroup
  \renewcommand\thefootnote{}\footnote{#1}%
  \addtocounter{footnote}{-1}%
  \endgroup
}
\journal{}
\begin{document}
\centerline{}
\begin{frontmatter}


\selectlanguage{english}
\title{The Scaling Site\\\vspace{.05in}
Le Site des Fr\' equences}


\selectlanguage{english}
\author[authorlabel1]{Alain Connes},
\ead{alain@connes.org}
\author[authorlabel2]{Caterina Consani\thanksref{label2}}
\thanks[label2]{The second author thanks the Coll\`ege de France for hospitality and financial support.}
\ead{kc@math.jhu.edu}

\address[authorlabel1]{Coll\`ege de France,
3 rue d'Ulm, Paris F-75005 France;
I.H.E.S. and Ohio State University}
\address[authorlabel2]{The Johns Hopkins
University, Baltimore, MD 21218 USA}


\begin{abstract}
\selectlanguage{english}

\vspace*{-.4cm}

We investigate the semi-ringed topos obtained  from the arithmetic site $\aarith$ of \cite{CC,CC1}, by extension of scalars from the smallest Boolean semifield $\B$ to the tropical semifield $\rmax$. The obtained  site $\rnt$ is the semi-direct product of the Euclidean half-line and the monoid $\nt$ of positive integers acting by multiplication. Its points are the same as the points $\aarith(\rmax)$  of $\aarith$ over $\rmax$ and form the quotient of the ad\`ele class space of $\Q$ by the action of the maximal compact subgroup $\hat{\Z}^*$ of the id\`ele class group. The structure sheaf  of the scaling topos endows it with a natural structure of tropical curve over the topos $\wnt$.  The restriction of this structure to the periodic orbits of the scaling flow gives, for each prime $p$, an analogue  $C_p$ of an elliptic curve whose Jacobian is $\Z/(p-1)\Z$. The Riemann-Roch formula holds on $C_p$ and involves real valued dimensions and real degrees for divisors.


\vskip 0.5\baselineskip

\selectlanguage{francais}

\noindent{\bf R\'esum\'e} \vskip 0.5\baselineskip \noindent

Le Site des Fr\' equences $\rnt$ est  obtenu \`a partir du site arithm\'etique $\aarith$ de \cite{CC,CC1} par extension des scalaires du semicorps bool\' een $\B$ au semicorps tropical $\rmax$. C'est le produit semi-direct de la demi-droite Euclidienne $[0,\infty)$  par l'action du semi-groupe $\nt$ des entiers positifs  par multiplication. Ses points sont les m\^emes que ceux du site arithm\' etique d\'efinis  sur $\rmax$ et forment le quotient de l'espace des classes d'ad\`eles de $\Q$ par l'action du sous-groupe compact maximal du groupe des classes d'id\`eles. 
Le faisceau structural du site des fr\'equences en fait une courbe tropicale dans le topos $\wnt$. La restriction de cette structure aux orbites p\' eriodiques donne, pour chaque nombre premier $p$, un analogue $C_p$ d'une courbe elliptique dont la Jacobienne  est  $\Z/(p-1)\Z$. La formule de Riemann-Roch pour $C_p$ fait apparaitre  des dimensions \`a valeurs r\'eelles et les  degr\'es des diviseurs  sont des nombres r\'eels.


\end{abstract}
\end{frontmatter}


\selectlanguage{english}
\vspace*{-1.4cm}
\section{Introduction}
\label{}
\blfootnote{{\it Keywords:}~Site Arithm\'etique, Site des fr\'equences, classes d'Ad\`eles, topos,  caract\'eristique 1.
Arithmetic Site, Scaling site, Ad\`ele class space, topos, characteristic 1.}

\vspace*{-.7cm}

This note describes the Scaling Site as the algebraic geometric space  obtained from the arithmetic site $\aarith$ of \cite{CC,CC1} by extension of scalars from the Boolean semifield $\B$ to the tropical semifield $\rmax$. The underlying  site $\rnt$ inherits, from its sheaf of regular functions, a natural structure of tropical curve allowing one to define the sheaf of rational functions and to   investigate an adequate version of the Riemann-Roch theorem in characteristic $1$. We test this structure by restricting it to the periodic orbits of the scaling flow, \ie to the points over the image of $\Spec\Z$ (\cf\cite{CC1}, \S 5.1). We find that for each prime $p$ the corresponding circle of length $\log p$ is endowed with a quasi-tropical structure which turns this orbit into the analogue $C_p=\R_+^*/p^\Z$ of a classical elliptic curve $\C^*/q^\Z$. In particular the notions of rational functions, divisors, etc are all meaningful. A new feature is that the degree of a divisor can now be any real number. We determine the Jacobian of  the curve $C_p$, \ie the quotient $J(C_p)$  of the group of divisors of degree $0$ by principal divisors and show  in Theorem \ref{thmjaccp}  that it is a cyclic group of order $p-1$. For each divisor $D$ on $C_p$ we define the corresponding Riemann-Roch problem with solution space $H^0(D):=H^0(C_p,\cO(D))$. We introduce the continuous dimension $\cdim(H^0(D))$ of this $\rma$-module using a limit of normalized topological dimensions and find that $\cdim(H^0(D))$ is a real number. Finally, in Theorem \ref{RRperiodic} we prove  that the  Riemann-Roch formula holds for $C_p$. The appearance of arbitrary positive real numbers as continuous dimensions in this formula is due to the density in $\R$ of the subgroup $H_p\subset \Q$ of fractions with denominators a power of $p$ and the fact that continuous dimensions are obtained as limits of normalized dimensions $p^{-n}\tdim(H^0(D)^{p^n})$. We view this outcome  as the analogue in characteristic $1$ of what happens for matroid  $C^*$-algebras and the type II normalized traces as in \cite{dix}.

\subsection{Notations} For any abelian ordered group $H$ we let $H_{\rm max}=H\cup \{-\infty\}$ be the semifield obtained from $H$ by applying the max-plus construction, \ie the addition is given by the max, and the multiplication by the addition in $H$. In particular $\R_{\rm max}$ is isomorphic to $\rmax$ by the exponential map (\cf \cite{Gaubert}). 
\vspace*{-.5cm}

\section{The scaling site}
\vspace*{-.3cm}

  The scaling site $\rnt$ is, as a  site, given by a small category $C$ endowed with a Grothendieck topology $J$. The objects of $C$ are the (possibly empty) bounded open intervals 
$\Omega\subset [0,\infty)$. The morphisms between two objects are defined by 
$
\Hom_C(\Omega,\Omega')=\{n\in \nt\mid n\Omega\subset \Omega'\},
$
 if $\Omega\neq \emptyset$ and by $\Hom_C(\emptyset,\Omega'):=\{*\}$ \ie the one point set, for any object of $C$. Thus the empty set is the initial object of $C$.
The category $C$ admits pullbacks. Indeed, let $\Omega_j\neq \emptyset$ ($j=1,2$) and consider two morphisms $\phi_j:\Omega_j\to \Omega$ given by integers $n_j\in \Hom_C(\Omega_j,\Omega)$. Let $n={\rm lcm}(n_j)$ be their lowest common multiple, write $n=a_jn_j$ and let $\Omega':=\{\lambda \in [0,\infty)\mid a_j\lambda \in \Omega_j, \ j=1,2\}$. If $\Omega'=\emptyset$ the initial object is the pullback. Otherwise this gives an object $\Omega'$ of $C$ and morphisms $a_j\in \Hom_C(\Omega',\Omega_j)$ such that $\phi_1\circ a_1=\phi_2\circ a_2$. One sees that $(\Omega',a_j)$ is the pullback of the pair $\phi_j:\Omega_j\to \Omega$. Since the category $C$ has pullbacks we can describe a Grothendieck topology $J$ on $C$ by providing a basis (\cf \cite{MM}, Definition III.2). 
\vsp

\begin{e-proposition}
$(i)$~For each object $\Omega$ of $C$, let $K(C)$ be the collection of all ordinary covers $\{\Omega_i\subset \Omega, i\in I\mid \cup \Omega_i=\Omega\}$ of $\Omega$. Then $K$ defines a Grothendieck topology $J$ on $C$.

$(ii)$~The category $\sh(C,J)$ of sheaves on $(C,J)$ is canonically isomorphic to the category of $\nt$-equivariant sheaves on $[0,\infty)$.
\end{e-proposition} 
\vsp
\begin{e-definition}
The {\em scaling site} $\rnt$ is the small category $C$ endowed with the Grothendieck topology $J$. The scaling topos is the category $\sh(C,J)$.
\end{e-definition}
\vsp
\vspace*{-.5cm}

\section{The points of the scaling topos} 
\vspace*{-.3cm}

We recall from \cite{CC1} that the space $\aarith(\rmax)$ of points of the arithmetic site $\aarith$ over $\rmax$  is the disjoint union of the following two spaces:

$(i)$~The points which are defined over $\B$:  they correspond to the  points of $\wnt$ and are in canonical bijection with the space $
\Q_+^\times\backslash\A^f/\hat\Z^*
$  of  ad\`ele classes whose archimedean component vanishes.

$(ii)$~The points of $\aarith(\rmax)\setminus \aarith(\B)$  are in canonical bijection with the space $\Q_+^\times\backslash\left(( \A^f/\hat\Z^*)\times \R_+^*\right)$ of  ad\`ele classes whose archimedean component does not vanish. Equivalently, these points correspond to the space $\mathfrak R $ of  rank one subgroups  of $\R$ through the map 
$$
(a,\lambda) \mapsto  \lambda H_a\qqq a\in \A_f/\hat\Z^*, \, \lambda \in  \R_+^*, \  \  H_a:=\{q\in\Q\mid qa\in \hat\Z \}.
$$
The next statement shows that the points of the scaling topos $\sh(C,J)$ are in canonical bijection with $\aarith(\rmax)$. We recall that the points of a topos of the form $\sh(C,J)$ are equivalently described as flat, continuous functors $F:C\to \Se$ (\cf \cite{MM} VII.6 Corollary 4). In our context, we define the support of such a functor as the complement of the union of the open intervals $I$ such that $F(I)=\emptyset$.
\vsp
\begin{theorem}\label{thmpointsrnt}
$(i)$~The category of points of the scaling topos  with support $\{0\}$ is the same as the category of points of $\wnt$.\newline
$(ii)$~The category of points of the scaling topos  with  support different from $\{0\}$ is canonically equivalent to the category of rank one subgroups of $\R$.

\end{theorem}
\vsp
The proof of the above theorem follows from the next four lemmas.
\vsp
\begin{lemma} \label{doublelem1} 
$(i)$~Let $H\subset \R$ be a rank one subgroup. Then  $F_H(V):=V\cap H_+$ defines a flat, continuous functor $F_H:C\to \Se$.

$(ii)$~The map $H\mapsto \ffp_H$ which associates to a rank one subgroup of $\R$ the point of $\sh(C,J)$ represented by the flat continuous functor $F_H$ is an injection of $\mathfrak R $ in the space of points of the scaling topos  up to isomorphism.
\end{lemma}
{\it\bf Proof.} $(i)$~One verifies that the category $\int_C F_H$ is filtering, \ie fulfills the three conditions of Definition VII.6.2 of \cite{MM}. The fact that $H$ is of rank $1$ yields the second filtering condition. The third condition is automatic since two morphisms $u,v\in \Hom(I,J)$ which fulfill $F_H(u)x=F_H(v)x$ for some $x\in F_H(I)$ are necessarily equal. Moreover the functor $F_H$ is continuous in the sense that its maps a covering  to an epimorphic family (\cf \cite{MM} Lemma VII. 5.3).

$(ii)$~Given a point $\lambda\in (0,\infty)$ let $\{V_j\}$ be open intervals forming a basis of neighborhoods of $\lambda$. Then one has
$
\cap F_H(V_j)\neq \emptyset \iff \lambda \in H
$
and this shows that one can recover the subgroup $H\subset \R$ from the continuous flat functor. This construction only depends upon the isomorphism class of the point $\ffp_H$.\qed
\vsp

The next lemma shows that the category of points of the scaling topos  with support $\{0\}$ is the same as the category of points of $\wnt$.\vsp

\begin{lemma}\label{case1} Let $F:C\to \Se$ be a flat continuous functor.
Assume  that  $F(V)=\emptyset$ when $0\notin V$. Then there exists a unique flat functor $X:\nt\to \Se$ such that $F(V)=X$ for any object $V$ of $C$ containing $0$.
\end{lemma}
{\it\bf Proof.}  Let  
$X:=\varprojlim F(J)$ where $J$ runs through the open intervals containing $0$. For any such interval $J$ one has a natural map $i_J:X\to F(J)$ and the continuity of $F$ shows that this map is bijective. Moreover $X$ inherits an action of $\nt$ which is uniquely specified by requiring \raggedbottom
$
i_J(n.x)=F(n)(i_{J/n}(x))\qqq x\in X, \  n\in \nt.
$
One checks that this construction makes sense independently of the choice of $J$, $0\in J$. Finally, the flatness of  
$F$ shows that the functor $\nt\to\Se$ obtained from the action of $\nt$ on $X$ is flat. \qed
\vsp

The next two lemmas show that the category of points of the scaling topos  with support $\neq \{0\}$ is equivalent to the category of rank one subgroups of $\R$.\vsp

\begin{lemma}\label{case2} Let $F:C\to \Se$ be a flat continuous functor.
Let $\lambda\in(0,\infty)$ and $F_\lambda:=\varprojlim_{\lambda\in J} F(J)$ be the co-stalk of $F$ at $\lambda$. Then there exists at most one element in the set $F_\lambda$ and for any bounded open interval $V\subset (0,\infty)$, $F(V)$ is the disjoint union $\cup_{\lambda \in V} F_\lambda$.
\end{lemma}
{\it\bf Proof.} 
 We first show that $F(V)=\cup_{\lambda \in V} F_\lambda$. Let $z\in F(V)$ then, by continuity of $F$, it follows that for any  
cover $V=\cup V_j$ one has $z\in F(V_j)$ for some $j$. One first writes  $V=\cup W_j$ with $W_j$ an increasing family of open intervals such that $\overline{ W_j}\subset W_{j+1}$. Then one gets an interval $W$, with $\overline{W}\subset V$ such that $z\in F(W)$. Using a family of covers $\cU_k$ of $W$ such that the maximal diameter of the open sets in $\cU_k$ tends to $0$, and is less than the Lebesgue number of $\cU_{k-1}$ one obtains a decreasing sequence of intervals $I_k\subset W$ such that $z\in F(I_k)$ for all $k$, and the unique element $\lambda\in\cap I_k\subset \overline{W}\subset V$ is such that $z\in F_\lambda$. Next we show that $F(V)=\cup_{\lambda \in V} F_\lambda$ is a disjoint union. Indeed if $U,U'$ are disjoint open intervals contained in $V$ one has $F(U)\cap F(U')=\emptyset$ inside  $F(V)$. Assume on the contrary that there exist $z\in F(U)$ and $z'\in F(U')$ such that  $F(\iota)z=F(\iota')z'$ where $\iota:U\to V$ and $\iota':U'\to V$ are the inclusions. By applying the flatness of $F$ let $W$ an object of $C$, $u\in F(W)$, $n,n'\in \nt$ be such that $nW\subset U$, $n'W\subset U'$ and $F(j)(u)=z$ for $j=(n:W\to U)$, $F(j')(u)=z'$ for $j'=(n':W\to U')$. The two morphisms $\iota\circ j$ and $\iota'\circ j'$, $W\to V$ fulfill $F(\iota\circ j)u=F(\iota'\circ j')(u)$. Then the third property of a flat functor shows that there exists a morphism in $C$ which equalizes $\iota\circ j$ with $\iota'\circ j'$. Thus $n=n'$ and one gets a contradiction since $U\cap U'=\emptyset$. 
 Next we  show that $F_\lambda$ contains at most one element.
Let $z,z'\in F_\lambda$. Then again by flatness of  $F$, there exist for any open interval $I\ni\lambda$ an object $W$ of $C$, an element $u\in F(W)$ and integers $n,n'\in \nt$ such that $nW\subset I$, $n'W\subset I$ and $F(j)(u)=z$ for $j=(n:W\to I)$ and $F(j')(u)=z'$ for $j'=(n':W\to I)$. Let $\mu\in W$ be the unique element  such that $u\in F_\mu$. Then one has $F(j)u\in F_{n\mu}$ and  by uniqueness one gets $n\mu=\lambda$. Similarly $n'\mu=\lambda$ and, since $\lambda\neq 0$, one has $j=j'$ and $z'=z$ so that $F_\lambda$ contains at most one element.\qed
\vsp

\begin{lemma}\label{case2bis} Let $F:C\to \Se$ be a flat continuous functor.
Assume that $F(V)\neq \emptyset$ for some open interval $V$ not containing $0$. Then the set $H^+_F:=\{\lambda\in(0,\infty)\mid F_\lambda\neq \emptyset\}$ is the positive part of a rank one subgroup $H_F$ of $\R$.
\end{lemma}
{\it\bf Proof.} By Lemma \ref{case2}, the subset $H^+_F$ is non-empty. For each $n\in\nt$ the multiplication by $n$ maps $H^+_F$ to itself using the morphism in $C$ given by $n:I\to nI$ in a small neighborhood of $\lambda$ with $F_\lambda\neq \emptyset$. Moreover the flatness of $F$ shows that given two elements $\lambda,\lambda'\in H^+_F$, there exists $\mu\in H^+_F$ and $n,n'\in\nt$ such that $n\mu=\lambda$ and $n'\mu=\lambda'$. It follows that $H^+_F$ is an increasing union of subsets of the form $h_k\Z\cap (0,\infty)$, $h_k>0$, and one gets $H^+_F=H_F\cap (0,\infty)$ where $H_F$ is the increasing union of the subgroups $h_k\Z$. \qed
\vspace*{-.5cm}
 
\section{The structure sheaf $\cO=\zmax\hat\otimes_\B\rmax$ and its stalks} 

\vspace*{-.3cm}

The Legendre transform allows one to describe the reduced semiring $\zmax\hat\otimes_\B\rmax$  involved in the extension of scalars of the arithmetic site $\aarith$ from $\B$ to $\rmax$ in terms of 
$\rma$-valued functions on $[0,\infty)$ which are convex, piecewise affine functions with integral slopes.
We first discuss an analogous result that holds when $\zmax$ is replaced by  the semiring $H_{\rm max}$ associated by the max-plus construction to a rank one subgroup $H\subset \R$.

\vspace*{-.3cm}
\subsection{The Legendre transform}
\vspace*{-.3cm}

Let us fix a rank one subgroup $H\subset \R$ and consider the tensor product 
$H_{\rm max}\otimes_\B\rma$ and the associated \mc semiring $R=H_{\rm max}\hat \otimes_\B\rma$ 
whose elements are viewed as Newton polygons with vertices  pairs $(x,y)\in H\times \R$ (\!\!\cite{CC1}). Let $Q=H_+\times \R_+$. Any element of $R$ is given by the convex hull $N$ in $\R^2$ of the union of finitely many quadrants $(x_j,y_j)-Q$. This convex hull $N$ is the intersection of half planes $P\subset \R^2$  of the form
$
P_{\lambda,u}:=\{(x,y)\mid \lambda x+y\leq u\}, \  P^v:=\{(x,y)\mid x\leq v\}
$, where $\lambda \in \R_+$ and $u,v\in\R$.
This description shows that $N$ is uniquely determined by the function 
$
\ell_N(\lambda):=\min \{u\in\R\mid N\subset P_{\lambda,u}\}
$
and that this function is given in terms of the finitely many vertices $(x_j,y_j)$ of the Newton polygon $N$ by the formula 
\begin{equation}\label{hdefn}
\ell_N(\lambda)=\max_j \lambda x_j+y_j.
\end{equation}

\begin{e-proposition} \label{legendrelem} Let $H\subset \R$ be a subgroup of rank one.
The map $N\mapsto \ell_N$ is an isomorphism of the \mc semiring $R=H_{\rm max}\hat\otimes_\B\rma$ with the semiring $\cR(H)$ of convex, piecewise affine continuous functions on 
$[0,\infty)$ with slopes in $H\subset \R$ and only finitely many singularities. The operations are the pointwise 
operations of $\rma$-valued functions.
\end{e-proposition}

\vspace*{-.3cm}
\subsection{The stalks of $\cO$}
\vspace*{-.3cm}

Proposition  \ref{legendrelem} gives the relation between the reduced semiring $\zmax\hat\otimes_\B\rmax$  involved in the extension of scalars of the arithmetic site from $\B$ to $\rmax$, and the semiring $\cR(\Z)$.
The structure sheaf $\cO$ of $\rnt$ is defined by localizing the semiring $\cR(\Z)$. The sections $\xi \in \cO(\Omega)$ on an open set $\Omega\subset [0,\infty)$ are convex, piecewise affine continuous functions on 
$\Omega$ with slopes in $\Z\subset \R$ and only finitely many singularities. The action of $\nt$ on $\cO$  is given by the morphisms 
\begin{equation}\label{gammar}
\gamma_n:\cO(\Omega)\to \cO(\frac 1n\Omega), \  \   \gamma_n(\xi)(\lambda):=\xi(n\lambda)\qqq \lambda \in [0,\infty), \  n\in \nt.
\end{equation}
For $\xi(\lambda)=\max\{\lambda h_j+s_j\}$ as in \eqref{hdefn} one has $\xi(n\lambda)=\max\{\lambda n h_j+s_j\}$ so that $\gamma_n(\xi)\in \cO(\frac 1n\Omega)$. Note that these maps are not invertible. \vsp

\begin{figure}
\centering
\begin{subfigure}{.5\textwidth}
  \centering
  \includegraphics[scale=0.3]{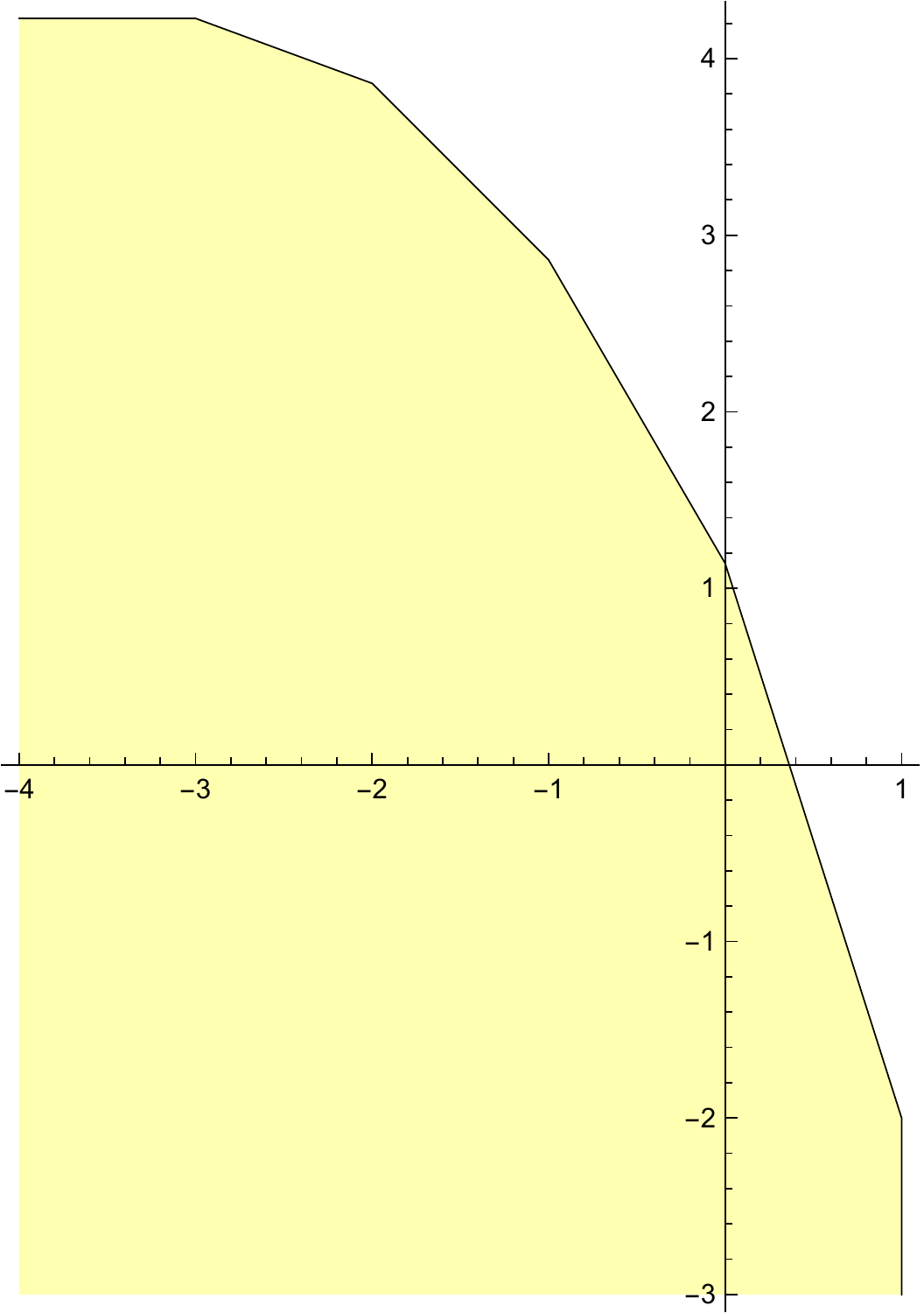}
  \caption{Element $C$ of $\cR(\Z)$}
  \label{elem}
\end{subfigure}%
\begin{subfigure}{.5\textwidth}
  \centering
  \includegraphics[scale=0.5]{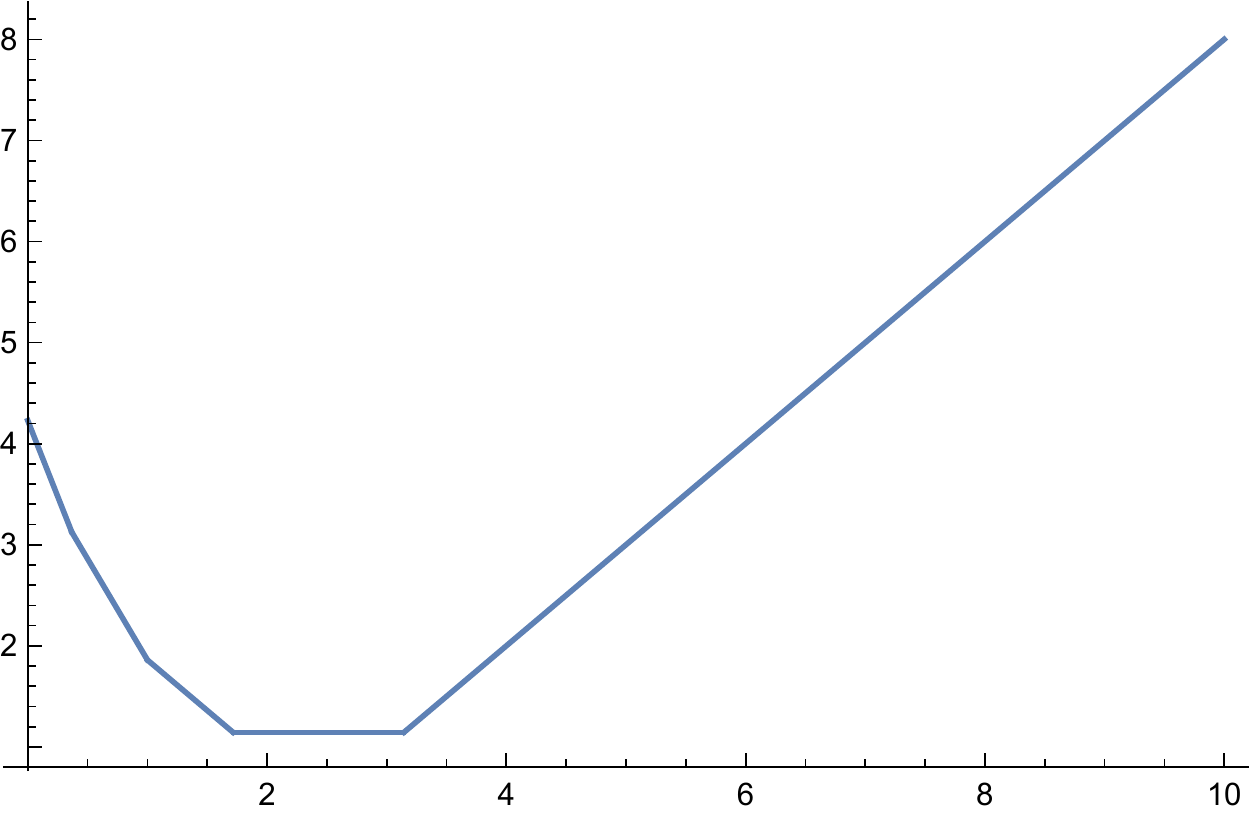}
  \caption{The Legendre transform $\ell_C(\lambda)$}
  \label{lege}
\end{subfigure}
\end{figure}

\begin{theorem} \label{structure2} $(i)$~Let  $H\subset \R$ be a rank one subgroup of $\R$ and $\ffp_H$ be the  associated point of the scaling topos. The stalk of the structure sheaf $\cO$ at $\ffp_H$  is the semiring $\cO_H$ of germs  of $\rma$-valued, piecewise affine, convex continuous functions with slope in $H$.

$(ii)$~Let $H$ be an abstract rank one ordered group and $\ffp_H^o$ the point  of the scaling topos with support $\{0\}$ associated to $H$. The stalk of the structure sheaf $\cO$ at $\ffp_H^o$ is the semiring $\cZ_H=(\R\times H)_{\rm max}$ associated by the max-plus construction to the totally ordered group $\R\times H$ endowed with the lexicographic order.
\end{theorem}
{\it\bf Proof.} 
$(i)$~To evaluate the stalk of the structure sheaf $\cO$ at the point $\ffp_H$ we use the description 
$\ffp_H=\varinjlim y_{I_j}$  as a filtered colimit of the Yoneda functors $y_{I_j}(V):=\Hom(I_j,V)$ where the elements $h_j\in H$, the objects $I_j$ and integers $n_j\in \nt$ fulfill the conditions ($\vert I\vert$ denotes the diameter of the interval $I$)
\begin{equation*}
H=\cup h_j\Z, \  n_jh_{j+1}= h_j,   \, h_j\in I_j, 
\,  n_j\bar I_{j+1}\subset I_j, \ \forall j\geq 1, \  \ 
\,\lim_{k\to \infty}(\prod_1^{k-1} n_i)\, \vert I_k\vert = 0
\end{equation*}
 The stalk of $\cO$ at the point $\ffp_H$ is 
$
\cO_{\ffp_H}=\varinjlim \cO(I_j)
$.
We define a map $\rho:\cO_{\ffp_H}\to \cR_H$ by associating to $(j,f)$, $f\in \cO(I_j)$ the germ of the function $\lambda\mapsto f(\lambda h_j)$ at $\lambda=1$. This function is defined in the neighborhood $\{\lambda\mid h_j\lambda \in I_j\}$ of $\lambda=1$. It is a piecewise affine convex continuous function with slopes in $h_j\Z\subset H$. Thus its germ  at $\lambda=1$  is an element  $\rho(j,f)\in \cR_H$. This construction is compatible with the colimit  in the sense  that $\rho(j,f)=\rho(j+1,\gamma_{n_j}(f))$ where $\gamma_n$ is defined in \eqref{gammar}. Indeed one has
$
\gamma_{n_j}(f)(\lambda)=f(n_j\lambda)\qqq \lambda \in [0,\infty)
$.
Thus, using $n_j h_{j+1}=h_j$, one obtains
$$
\rho(j+1,\gamma_{n_j}(f))(\lambda)=\gamma_{n_j}(f)(\lambda h_{j+1})=f(n_j\lambda h_{j+1})=f(\lambda h_j)=\rho(j,f)(\lambda)
$$
One derives an isomorphism of semirings $\rho:\cO_{\ffp_H}\to \cR_H$.

$(ii)$~Let $H:=\varinjlim \Z$ where we use the $n_j$'s to organize the inductive system. We denote by $\iota(j,k)$ the element of $H$ associated to the image of $k\in \Z$ by the canonical map $\Z\to H$ associated to the $j$-th copy of $\Z$ in the colimit. By construction one has the equality
$
\iota(j,k)=\iota(j+1, n_j\, k)\qqq j, \ k\in \Z
$.
Then the stalk of the structure sheaf $\cO$ at the point $\ffp_H^o$ is the colimit $\cO_{\ffp_H^o}=\varinjlim \cO(I_j)$. We define a map $\delta:\cO_{\ffp_H^o}\to H$ as follows. We associate to $(j,f)$, $f\in \cO(I_j)$, the element 
$\delta(j,f):=\iota(j,k)$ where $k=f'(0)\in \Z$ is the derivative of $f$ at $0\in I_j$. One then has 
$$
\delta(j+1,\gamma_{n_j}(f))=\iota(j+1,\gamma_{n_j}(f)'(0)=\iota(j+1,n_jf'(0))=\iota(j,f'(0))=\delta(j,f).
$$
This shows that $\delta:\cO_{\ffp_H^o}\to H$ is well defined. Similarly $\alpha(j,f):=f(0)$ defines a map $\alpha:\cO_{\ffp_H^o}\to \rma$ and the pair $\rho=(\alpha,\delta)$ determines a map $\cO_{\ffp_H^o}\to \cZ_H$ which is both injective and surjective. One checks that this map is an isomorphism of $\cO_{\ffp_H^o}$ for the semiring structure whose multiplication is given by $(x,h)\bullet (x',h')=(x+x',h+h')$ and addition is defined as
$$
(x,h)\vee (x',h'):=\begin{cases} (x,h)\ \text{if}\ x>x'\\(x',h')\ \text{if}\ x'>x\\
(x,h \vee h') \ \text{if}\ x=x'\end{cases}
$$  
\qed

 The germs at $\lambda=1$ of $\rma$-valued, piecewise affine, convex continuous functions $f(\lambda)$ with slopes in $H$ are characterized by a triple $(x,h_+,h_-)$,  such that $f(1\pm \epsilon)=x
\pm h_\pm \epsilon$ for $\epsilon \geq 0$ small enough. Here, one has $x\in \R$, $h_\pm\in H$, $h_+\geq h_-$. The only additional  element of the semiring $\cR_H$ corresponds to the germ of the constant function $-\infty$. This function is the zero element of the semiring. The algebraic rules  for non-zero elements in $\cR_H$ are as follows. The addition $\vee$ in $\cR_H$ is given by the max of the two germs:
$$
(x,h_+,h_-)\vee (x',h'_+,h'_-):=\begin{cases} (x,h_+,h_-)\ \text{if}\ x>x'\\(x',h'_+,h'_-)\ \text{if}\ x'>x\\
(x,h_+\vee h'_+, h_-\wedge h'_-) \ \text{if}\ x=x'\end{cases}
$$
The product in $\cR_H$ is given by the sum of the two germs 
$
(x,h_+,h_-)\bullet (x',h'_+,h'_-):=(x+x',h_++h'_+,h_-+h'_-).
$
When $f$ is viewed as a locally defined map $H\mapsto f(H)\in\rma$ from rank one subgroups of $\R$ to $\rma$, the associated germ $(x,h_+,h_-)$ of $f$ at $H$ is given by $x=f(H)$, $h_
\pm=\lim_{\epsilon\to 0\pm}(f((1+\epsilon) H)-f(H))/\epsilon$. 
\vsp

We shall denote by $\scal1$ the semi-ringed topos $(\rnt,\cO)$. We  view it as a relative topos over $\rmax$ in the sense that the structure semirings are over $\rmax$. Likewise for the arithmetic site the structure sheaf has no non-constant global sections.
\vspace*{-.5cm}

\section{The points of  $\scal1$  over $\rmax$}
\vspace*{-.3cm}

The next  Theorem states that extension of scalars from $\aarith$ to $\scal1$ does not affect the points over $\rmax$.\vsp

\begin{theorem} \label{structure3} The canonical projection from points of $\scal1$ defined over $\rmax$ 
to points of the scaling topos  is bijective.
\end{theorem}
\vsp
The proof of this theorem follows from Theorem \ref{structure2} and the following lemma.
\vsp
\begin{lemma}\label{homtormax} $(i)$~The map $(x,h_+,h_-)\mapsto x$ is the only element of $\Hom_{\rma}(\cR_H,\rma)$.
\newline
$(ii)$~The map $(x,h_+)\mapsto x$ is the only element of $\Hom_{\rma}(\cZ_H,\rma)$.
\end{lemma}

\vspace*{-.5cm}

\section{The real valued Riemann-Roch Theorem on periodic orbits}

\vspace*{-.3cm}
To realize the notion of rational functions in our context we proceed as in the definition of Cartier divisors and consider the sheaf obtained from the structure sheaf $\cO$ by passing to the semifield of fractions.

\vsp

\begin{e-proposition}\label{tensring} 
For any object $\Omega$ of $C$ the semiring $\cO(\Omega)$ is \mc and the canonical morphism to its semifield of fractions $\cK(\Omega)$ is the inclusion of convex, piecewise affine, continuous functions among continuous, piecewise affine functions, endowed with the two operations of max and plus.
\end{e-proposition}\vsp
The natural action of $\nt$ on  $\cK$ defines a sheaf of semifields in the scaling topos. One determines its stalks in the same way as for the structure sheaf $\cO$. The local convexity no longer holds, \ie the difference $h_+-h_-\in H\subset \R$ is no longer required to be positive.
\vsp
\begin{e-definition}\label{site2} Let $\ffp_H$ be the point of the scaling topos associated to the rank one subgroup $H\subset \R$ and let $f$ be an element of the stalk of $\cK$ at $\ffp_H$. The order of $f$ at $H$ is defined as ${\rm Order}(f)=h_+-h_-\in H\subset \R$ where $h_
\pm=\lim_{\epsilon\to 0\pm}(f((1+\epsilon) H)-f(H))/\epsilon$.
\end{e-definition}
\vsp
Let $p$ be a prime and consider the subspace $C_p$ of points of  $\rnt$ corresponding to subgroups $H\subset \R$ which are abstractly isomorphic to the subgroup $H_p\subset \Q$ of fractions with denominator a power of $p$. 
\vsp
\begin{lemma}\label{periodp} The map $\R_+^*\to C_p$,  $\lambda\mapsto\lambda H_p$ induces a topological isomorphism $\eta_p: \R_+^*/p^\Z\to C_p$. The pullback by $\eta_p$ of the structure sheaf $\cO$ is the sheaf $\cO_p$ on $\R_+^*/p^\Z$ of piecewise affine, continuous convex functions, with slopes in $H_p$.
\end{lemma}
\vsp
We use $\eta_p$ to view functions on $C_p$ as functions of $\lambda\in \R_+^*/p^\Z$. Note that at $H=\lambda H_p$ one has  $h_
\pm=\lim_{\epsilon\to 0\pm}(f((1+\epsilon) H)-f(H))/\epsilon=\lambda f'^\pm(\lambda)$ where $f'^\pm(\lambda)$ are the directional derivatives, and that the condition $h_\pm\in H$ means that $\lambda f'^\pm(\lambda)\in \lambda H_p$ \ie $f'^\pm(\lambda)\in  H_p$. We now apply the notion of order as in Definition \ref{site2} to the global sections of the sheaf of quotients of the sheaf of semirings $\cO_p$.
\vsp
\begin{lemma}\label{periodpbis} $(i)$~The sheaf of quotients $\cK_p$ of the sheaf of semirings $\cO_p$ is the sheaf  (on $\R_+^*/p^\Z$) of piecewise affine, continuous  functions with slopes in $H_p$,  endowed with the two operations of max and plus. \newline
$(ii)$~The sheaf $\cK_p$ admits global sections and for any $f\in H^0(\R_+^*/p^\Z,\cK_p)$ one has: 
\begin{equation*}
\sum_{\R_+^*/p^\Z} {\rm Order}(f)(\lambda)=0.
\end{equation*}
\end{lemma}
A divisor  $D$ is a section $H\mapsto D(H)\in H\subset \R$, vanishing except on a finite subset, of the projection on the base from the total space of the bundle formed by pairs $(H,h)$ where $H\subset \R$ is a subgroup abstractly isomorphic to the subgroup $H_p\subset \Q$ and where $h\in H$. The degree of a divisor $D$ is the finite sum $\deg(D):=\sum_{H\in C_p} D(H)\in\R$. 
Next, we define an invariant of divisors with values in the group $H_p/(p-1)H_p\simeq \Z/(p-1)\Z$. Note that given $H\in C_p$, the elements $\lambda\in \R_+^*$ such that $H=\lambda H_p$ determine maps $\lambda^{-1}:H\to H_p$ differing from each other by multiplication by a power of $p$, thus  the corresponding map $\chi:H\to H_p/(p-1)H_p\simeq \Z/(p-1)\Z$ is canonical. For any divisor $D$ on $C_p$ we define
\begin{equation*}
\chi(D):=\sum_{H\in C_p} \chi(D(H))\in \Z/(p-1)\Z.
\end{equation*}
Then, we obtain the following\vsp
\begin{theorem}\label{thmjaccp} The map $\chi: H\to H_p/(p-1)H_p$ vanishes on principal divisors and it induces an isomorphism of groups $\chi:J(C_p)\to \Z/(p-1)\Z$ of the quotient  $J(C_p)=\div(C_p)^0/\cP$ of the group of divisors of degree $0$ by the subgroup $\cP$ of principal divisors.	
\end{theorem}
\vsp
Since the group law on divisors is given by pointwise addition of sections, both the maps $\deg:\div(C_p)\to \R$ and $\chi:\div(C_p)\to \Z/(p-1)\Z$ are group homomorphisms and the pair $(\deg,\chi)$ provides an isomorphism of groups 
\begin{equation}\label{jacmapbis}
(\deg,\chi):\div(C_p)/\cP\to \R\times (\Z/(p-1)\Z).
\end{equation}
 Given a divisor $D\in \div(C_p)$ one defines the following module over $\rmax$:
$$
H^0(D)=\Gamma(C_p,\cO(D)) =\{f\in \cK(C_p)\mid D+(f)\geq 0\}.
$$
\begin{e-definition}\label{pnorm} Let $f\in \Gamma(C_p,\cK_p)$.  	
 One sets  	
\begin{equation}\label{defnnp}
\Vert f\Vert_p :=\max \{\vert h(\lambda)\vert_p/\lambda\mid \lambda \in C_p\}\end{equation}
where $h(\lambda)\in H_p$ is the slope\footnote{at a point of discontinuity of the slopes one takes the max of the two values $\vert h_\pm(\lambda)\vert_p/\lambda$ in \eqref{defnnp}} of $f$ at $\lambda$.
\end{e-definition}
\vsp
Let $D\in \div(C_p)$ be a divisor. We introduce the following  increasing filtration of $H^0(D)$ by $\rma$-submodules:
\begin{equation*}
H^0(D)^\rho:=\{f\in H^0(D)\mid \Vert f\Vert_p\leq \rho\}.
\end{equation*}
\vsp
We denote by $\tdim(\cE)$ the topological covering dimension   of an $\rma$-module $\cE$ (\cf \cite{Pears}) and define
\begin{equation}\label{rr1}
\cdim(H^0(D)):=\lim_{n\to \infty} p^{-n}\tdim(H^0(D)^{p^n}).
\end{equation}
One shows that the above limit exists: indeed, one has the following
\vsp
\begin{theorem}\label{RRperiodic}
$(i)$~Let $D\in \div(C_p)$ be a divisor with $\deg(D)\geq 0$. Then the limit in \eqref{rr1} converges and one has  
$\cdim(H^0(D))=\deg(D)$.\newline
$(ii)$~The following Riemann-Roch formula holds
\begin{equation*}
\cdim(H^0(D))-\cdim(H^0(-D))=\deg(D),\qquad \forall D\in \div(C_p).
\end{equation*}	
\end{theorem}
\vsp
One can compare the above Riemann-Roch theorem with the tropical Riemann-Roch theorem of \cite{BN,GK,MZ} and its variants. More precisely, let us consider an elliptic tropical curve $C$, given by a circle of length $L$. In this case, the structure of the group $\div(C)/\cP$ of divisor classes is inserted into an exact sequence of the  form
$
0\to \R/L\Z\to \div(C)/\cP\stackrel{\deg}{\to} \Z\to 0
$ (\cf  \cite{MZ}\!\! ).
This sequence is very different from the split exact sequence associated to $C_p$ and
deduced from \eqref{jacmapbis}, \ie 
$
0\to \Z/(p-1)\Z\to \div(C_p)/\cP\stackrel{\deg}{\to} \R\to 0
$.
The reason for this difference is due to the nature of the structure sheaf of $C_p$, when this sheaf is written in terms of the variable $u=\log\lambda$. This choice is dictated by the requirement that the periodicity condition $f(px)=f(x)$ becomes translation invariance by $\log p$. The condition for $f$ of being piecewise affine in the parameter $\lambda$  is expressed in the variable $u$ by the piecewise vanishing of $\Delta'f$, where $\Delta'$ is the elliptic translation invariant operator 
$
\Delta'(f):=\left(\frac{\partial}{\partial u}\right)^2\, f - \frac{\partial}{\partial u}\, f
$.
In terms of the variable $\lambda=e^u$, this operator takes the form $\lambda^2\left(\frac{\partial}{\partial \lambda}\right)^2$ and this fact explains why the structure sheaf of $C_p$ is considered as tropical (in terms of the variable $\lambda$).









\end{document}